\documentclass{amsart} 
\usepackage{amsthm,amssymb,amsmath,dsfont,bigdelim}
\usepackage{tikz}
\usetikzlibrary{arrows,automata}
\usepackage{graphicx}
\usepackage{float}
\usepackage[arc,all]{xy}
\usepackage{longtable}
\newtheorem{thm}{Theorem}[section]
\newtheorem{cor}[thm]{Corollary}
\newtheorem{lem}[thm]{Lemma}
\newtheorem{prop}[thm]{Proposition}
\theoremstyle{definition}
\newtheorem{defn}[thm]{Definition}
\newtheorem{exa}[thm]{Example}

\newtheorem*{proof M}{\textbf{Proof of  Main Theorem}}
\theoremstyle{remark}

\numberwithin{equation}{section}
\newcommand{\diam}{\mathrm{diam}}
\newcommand{\girth}{\mathrm{girth}}
\begin{document}
\title[non-commuting graph]{On the  non-commuting graph associated to a finite-dimensional Lie algebra }%
\author[A. Chareh khah moghhaddam]{Akram Chareh khah moghaddam}%
\address{Department of Mathematics, Faculty of Mathematical Sciences, Ferdowsi University of Mashhad, Mashhad, Iran }
\email{charehkhahmoghaddam.akram@gmail.com}%
\author[A. Erfanian]{Ahmad Erfanian}
\email{erfanian@um.ac.ir }
\address{Department of Pure Mathematics and Center of Excellence in Analysis on Algebraic
Structures, Faculty of Mathematical Sciences, Ferdowsi University of Mashhad, Mashhad, Iran }
%\thanks{0}%
\author[A. Shamsaki]{Afsaneh Shamsaki}
\email{shamsaki.afsaneh@yahoo.com}
\address{Department of Mathematics, Faculty of Mathematical Sciences, Ferdowsi University of Mashhad, Mashhad, Iran }
\keywords{Non-commuting graph, Lie algebra}%
\subjclass{17B60 , 05C60, 05C10}
%\keywords{Schur multiplier, nilpotent Lie algebra}
%\date{0}%
%\dedicatory{0}%
%\commby{0}%
% ----------------------------------------------------------------
\begin{abstract}
In this paper, we define  the non-commuting graph  associated to a Lie algebra $ L $  and obtain some basic graph properties such as connectivity, diameter, girth, Hamiltonian and Eulerian.  Moreover,  planarity, outer planarity and isomorphism between two
such graphs are also discussed in the paper.
\end{abstract}
\maketitle
 \section{Introduction}
Research topics in algebraic structures and graph theory have led to several exciting questions and results.
Some papers have been written about assigning a graph to a group or ring and determining its algebraic properties. 
One can associate a graph to a group $G $ by various ways. For example, intersection graph, non-commuting graph and relative non-commuting graph are defined(see \cite{Ak2, A, Tol}). One of the associated graphs to a group $ G $ is the non-commuting graph, whose vertices are $ G\setminus Z(G) $ and two distinct vertices $ x $ and $ y $
are adjacent  if and only if $ xy\neq yx. $
Some results in group theory may be extended to Lie theory.
These conclusions  are not always similar and even they are different in many situations.
Recently, Lie algebras have  attracted some authors due to wide applications in other sciences, especially  physics. Understanding the structure of fundamental particles and the formulation of quantum theories in high-energy physics make it inevitable for physicists to be familiar with Lie algebras.
%Also, a close relation exists between group theory and Lie theory. So,
\\
 In this paper, a graph $G_L $ is associated to a Lie algebra $ L $   and the properties of  this kind of graph
are investigated.
 Let $ L $ be a Lie algebra and $ Z(L) $ be its centre. We defined  $ G_L $ as follows: \\
 Take $ L\setminus Z(L) $ the vertices of $ G_L $ and join two distinct vertices $ x $ and $ y $ whenever $ [x, y]\neq 0. $ We call $ G_L $
 the non-commuting graph of a Lie algebra $ L. $ All graphs are considered undirected and simple (with no loops or multiple edges). For every graph, we denote
 the set of vertices and the edges of  $ G_L $ by $ V(G_L) $ and $ E(G_L), $ respectively. If $ L $ is abelian, then $ G_L $ is null
 graph. Our main goal is to detect the properties of a Lie algebra $ L $  with its associated graph $G_L $ and vice verse.\\
 In  section $ 2, $ we study some properties of $G_L$ for a non-abelian Lie algebra. In addition, some results are obtained about connectivity, Hamiltonian, Eulerian properties, the girth and diameter of  $ G_L. $  In  section $ 3, $ we determine under what conditions of $L $ the graph $ G_L $ is  planar and  outerplanar. In the last section, we consider  isomorphism between two  non-commuting graphs and show that under some conditions two Lie algebras have equal elements.\\
 %its relation with the degree of a vertex are obtained.\\
 In the rest of this section, we recall some concepts of  graphs and Lie algebras  that are used in the next sections.\\
Let $ \Gamma $ be an arbitrary graph.  The degree  of a vertex $ v\in V(\Gamma)$ denoted by  $ \deg(v), $  is the number of edges,  which are incident to $ v. $  The minimum and maximum degrees are be denoted by $ \delta(\Gamma) $ and $ \Delta(\Gamma), $ 
 respectively. Also,
 $ | V(\Gamma)| $  is called the size of $ \Gamma. $ If the degree of all vertices $ \Gamma$ is equal to $ | V(\Gamma)| - 1,$ then $\Gamma$ is a complete graph. This graph with $ n $ vertices is denoted by $ K_n. $ A graph $ \Gamma $ is regular if the degrees of all vertices of $ G $ are the same. 
 A path of length $n$ consist of vertices from $ v_0 $ to $ v_n $ and 
   a sequence of distinct  $ e_i $ such that $ e_i=\lbrace v_{i-1}, v_i \rbrace $  and $ 1\leq i \leq n. $ 
 A cycle is a path in a graph that starts and ends at the same vertex. A graph is connected when there is at least a path for every pair of vertices.
 If $ v $ and $ w $ are two vertices of $ \Gamma, $
 then $ d(v, w) $ is the length of the shortest path between $ v $ and $ w. $ The longest distance between  all pairs of vertices is
 called the diameter of $ \Gamma$ and denoted by $ \diam(\Gamma). $  The girth of $ \Gamma, $ $ \girth(\Gamma), $ is  the length of the shortest cycle in graph $ \Gamma.$ 
  In graph theory, a planar graph is a graph that can be drawn such that no the edges cross each other. The 
 following theorem is an important tool in proving some theorems related to the planar graphs in section $ 3. $
 \begin{thm}\cite[Corollary 9.5.3]{Bondy}\label{deg1}
 If $ \Gamma$ is a simple and  planar, then $ \delta(\Gamma)\leq 5. $ 
 \end{thm}   
If a graph can be embedded in the plane such that all vertices of the graph on the outer face is called an outerplanar graph. 
Every outerplanar graph is planar,
but the converse of fact is not true. For example, the complete graph with four vertices
is the smallest planar graph that is not outerplanar.
\begin{thm}\cite[Page 1083]{dege}\label{outer}
Every subgraph of an outerplanar contains a vertex with degree at most $ 2. $ 
\end{thm}
An Euler cycle is cycle that visits every edge of a graph exactly once. If graph $ G $ has an Euler cycle, then it is called Eulerain. 
The next theorem states a condition equivalent to an Eulerian graph. 
\begin{thm}\cite[Theorem 4.1]{Bondy}\label{eul}
A connected non-empty graph is  Eulerain if and only if it has no vertices of odd degree.
\end{thm}
A Hamiltonian cycle is a cycle that  visits each vertex  exactly once. Every graph that has a Hamiltonian cycle is called Hamiltonian graph. The following theorem is an important equment for detecting Hamiltonian graphs. 
\begin{thm}(Dirac's Theorem)\cite[Theorem 4.3]{Bondy}\label{Dir}
Let $ \Gamma $ be a simple graph of order $ n $ such that $ n\geq 3 $ and the degree of every vertex is at least $\frac{n}{2}. $ Then $ \Gamma $ is Hamiltonian.
\end{thm}
Recall that a bipartite graph is a graph if its vertex set can be divided into two non-empty and disjoint subsets
 $U$ and $V,$ that every vertex in $ U $ is adjacent to every vertex in $  V.$  
If $\vert U \vert = m $ and $\vert V \vert = n, $ then a bipartite graph is denoted $K_{m,n}.$
The next theorem is used  to detect a planar graph.
%shows that the non-commuting graph $G_L$ is not complete bipartite. 
\begin{thm}(Kuratowski's Theorem)\cite[Theorem 9.10]{Bondy}\label{Lo}
A  graph is planar if and only if it  contains no  subdivision of $ K_5 $ or $ K_{3,3}. $ 
\end{thm}
The neighbourhood of a vertex $ x\in V(G_L) $ is the set of all vertices which are adjacent  to $ x. $ A dominating set $ D $ for a graph $ G_L $
is a subset of $ V(G_L) $ such that every vertices of $ G_L $ is either $ D $ or  neighbourhood of $ D. $ The domination number $ \gamma (G_L) $
is the size of smallest dominating set of $ G_L. $\\
A Lie algebra $ L $  is a vector space over the field $ \mathbb{F} $ together with a bilinear map called the Lie bracket $ L\times L\rightarrow L $ with $ (x, y)\mapsto [x, y], $ satisfying the following axioms:
\begin{itemize}
\item[(i)]
$[x, x]=0\quad \text{for all} \quad x\in L,$
\item[(ii)]
$[x, [y, z]]+[z,[x, y]]+[y,[z, x]]=0 \quad \text{for all} \quad x, y, z\in L.$ 
\end{itemize}
%\begin{align*}
%&[x, x]=0\quad \text{for all} \quad x\in L, \cr
%&[x, [y, z]]+[z,[x, y]]+[y,[z, x]]=0 \quad \text{for all} \quad x, y, z\in L. 
%\end{align*} 
The Lie bracket $ [x, y] $ is often referred to as  the commutator of $ x $ and $ y. $ The second condition  is known as the Jacobi identity. Since 
the Lie bracket $ [-, -] $ is bilinear, we have $ 0=[x+y, x+y]=[x, x]+[x, y]+[y, x]+[y, y]. $ Hence the first condition  implies $ [x, y]=-[y, x]$ 
for all $ x, y\in L. $  For example, every vector space $ L $ with  $ [x, y]=0 $ for all $ x, y \in L $ is a Lie algebra. This  is called an abelian Lie algebra.  A subalgebra of $ L $ is a subspace  $ K $ such that is close under the Lie bracket. 
For instance, $ C_{L}(x)=\lbrace y\in L \vert [x, y]=0 \rbrace $  is the centralizer of element $ x $ in $ L $ and a subalgebra of $ L. $ 
An ideal of $ L $ is a subalgebra $ I $ with this condition $ [x, y] \in I, $ where $ x\in L $ and $ y\in L. $ 
Let $ I $  and $ J $ are two ideals of a Lie algebra $ L. $ It is obvious that 
 $[I, J]=\langle [x, y]\mid x\in I, y\in J \rangle $ is an ideal of  $ L. $ In particular, the derived subalgebra $ L^{2} $ is the ideal generates by all $ [x, y] $ such that $ x, y\in L. $ Another important ideal  is 
the centre of a Lie algebra $ L $  and defined as $ Z(L)=\lbrace x\in L \mid [x, y]=0~\text{for~~all}~ y\in L\rbrace. $
If let $ \dim L=n $ and $ \lbrace x_1, \dots, x_n\rbrace $ is a basis of $ L, $ then every element $ x $ of $ L $ of the form
$ x=\alpha_1 x_1+\dots +\alpha_n x_n, $ where $ \alpha_1, \dots, \alpha_n \in \mathbb{F}_q$ and the number of elements such as is equal to $ q^{n}.$ 
%\begin{equation*}
%Z(L)=\lbrace x\in L \mid [x, y]=0~\text{for~~all}~ y\in L\rbrace.
%\end{equation*}
All notations and terminologies are standard and reader can refer to \cite{wil, Bondy}. Moreover, from now on
 all Lie algebras are considered finite-dimensional over the field $\mathbb{F}_q$ unless otherwise stated. 
\section{Basic results}  
 In this section, we repeat the definition of  non-commuting graph of Lie algebra $ L $ over the field $ \mathbb{F}_q $ and then investigate 
 about  connectivity, diameter, girth, Hamiltonian, and domination number of this graph. First, we state the definition of the new
 graph as the following.
 \begin{defn}
 Let $ L $ be a Lie algebra and $ Z(L) $ be the centre of $ L. $ The non-commuting graph of $ L, $ denoted by $ G_L, $  is an undirected simple graph whose vertices are all elements in $L\setminus Z(L)$ and two distinct vertices $x$ and $y$ are adjacent if and only if $[x,y]\neq 0.$ 
 \end{defn}
\begin{lem}\label{deg}
Let $L$ be a Lie algebra and $G_L$ be the non-commuting  graph associated to $L.$ Then $\deg(x) = | L | - | C_L (x) |$ for every $x \in V(G_L).$
\end{lem}
\begin{proof}
Since $\deg(x) =| \lbrace y \in V(G_L) \vert [x,y] \neq 0\rbrace |  $ and $V(G_L)=L\setminus Z(L),$ we have 
\begin{align*}
\deg(x) &=  | V(G_L) \setminus C_L(x) |= ( | L  \setminus  Z(L) | ) - | C_L(x) |
\cr
&=| L |-  | C_L(x) \cup Z(L) |
= | L |- | C_L(x) |.
\end{align*}
Hence $\deg(x) =  | L | - | C_L(x) | $ for all $x \in V(G_L). $
\end{proof}

\begin{lem}\label{single}
Let $ L $ be a Lie algebra. Then  $ G_L $ does not  have any isolated vertex. 
\end{lem}
\begin{proof}
By contrary, let $ x $ be an isolated vertex of $ G_L. $ Then $ \deg(x)=0. $ We know that $\deg(x) =  | L | - | C_L(x) |, $ by Lemma 
\ref{deg}. Therefore  $  L = C_L(x)  $ and so $ x\in Z(L). $ It is  a contradiction with $ x\in V(G_L). $
\end{proof}
 In the following proposition, we 
show that the non-commuting graph $ G_L $ is connected.
\begin{prop}
Let $ L $ be a Lie algebra. Then $ G_L $ is connected. 
\end{prop}
\begin{proof}
 Let $ x $ and $ y $ be two arbitrary  distinct vertices of $ G_L. $ If they are  adjacent, then $ G_L $ is connected.
Suppose that   $x $ and $ y $ are  not adjacent. Then $ [x, y]=0. $ Since $ x, y \in V(G_L), $ there exist 
$ x_1, y_1 \in V(G_L) $ such that $ [x, x_1]\neq 0 $ and   $ [y, y_1]\neq 0 $ by Lemma \ref{single}. If $ [x,y_1] \neq 0 $ or $ [y, x_1]\neq 0, $   then $ x\sim y_1\sim y $ or $ x\sim x_1\sim y. $ Therefore
$G_L$ is connected. Now, let $ [x, y_1]= [y, x_1] =0 $ Then $ [x, x_1+y_1]=[x, x_1]+[x, y_1]\neq 0 $ and 
 $  [y, x_1+y_1]=[y, x_1]+[y, y_1]\neq 0. $ So, $ x_{1} + y_{1} \in V(G_L) $ and $ x\sim x_1+y_1\sim y $ is a path. Hence $ G_L $ is connected.
\end{proof}
The next proposition states the girth of non-commuting graph $ G_L $ is equal to $ 3. $
\begin{prop}
Let $ L $ be a Lie algebra. Then $ \girth (G_L)=3. $
\end{prop}
\begin{proof}
Assume that $ x $ is an arbitrary vertex of $ G_L. $ Then there exists $ y \in V(G_L) $ such that $ x $ and $ y $ are adjacent by Lemma \ref{single}. Since
$ [x, y]\neq 0, $ we have 
$[x+y, y]\neq 0$ and $[x+y, x]\neq 0.$
So, $ x+y \in V(G_L) $ and  $ x+y, $ $ x $  and $ y $ are adjacent. Hence there is a cycle of length $ 3 $ among
the vertices $ \lbrace x, y, x+y \rbrace $ and the girth of $ G_L $ is equal to $ 3. $
\end{proof}
In the next proposition, we obtain an upper bound for the diameter of non-commuting graph $ G_L. $
\begin{prop}
Let $ L $ be a Lie algebra. Then $ \diam( G_L)\leq 2. $ 
\end{prop}
\begin{proof}
Assume that $ x, y \in V(G_L). $ We show that $ d(x, y)\leq 2. $
If $ x$ and $ y $ are adjacent, then $ d(x, y)=1< 2. $ Suppose that $ x$ and $ y $ are not adjacent. Thus there are $ x_1, y_1 \in V(G_L) $
such that $ [x, x_1]\neq 0 $ and $ [y, y_1]\neq 0, $ by Lemma \ref{single}. Now, we have the following cases:\\
\textbf{Case 1.} $ x $ is adjacent to $ y_1. $ In this case,  $ d(x, y)=d(x, y_1)+d(y_1, y)=2\leq 2. $ 
\\
\textbf{Case 2.} 
$ y $ is adjacent to $ x_1. $ By a similar way, one can see $ d(x, y)\leq 2. $
\\
\textbf{Case 3.} 
$ x $ and $ y $ are not adjacent to $ y_1 $ and $ x_1, $ respectively. So, $ [x, y_1]=[y, x_1]=0. $  Since 
$[x, x_1+y_1] \neq  0$ and
$[y, x_1+y_1] \neq  0,$
we have $ x_1+y_1\in V(G_L) $ and 
\begin{equation*}
d(x,y)=d(x, x_1+y_1)+d(x_1+y_1, y)=2\leq 2. 
\end{equation*}
  Hence $ \diam(G_L)\leq 2, $ by the previous cases.
\end{proof}
In the following, we show that under what condition $ \diam(G_L)=2. $ First of all, we need to prove the next  theorem.
\begin{thm}\label{thm2.8}
Let $ L $ be a Lie algebra  over the field $ \mathbb{F}_q $ and $ G_L $ be complete. Then $ | Z(L) |=1 $ and $ q=2. $  
\end{thm}
\begin{proof}
Since $ G_L $ is complete, then  $ \deg (x)= |  V(G_L)|-1$ and so $\deg(x)=| L |-| Z(L)|-1 $ for all $ x\in V(G_L). $ On the other hand,  we know that
$ \deg (x)= | L |-| C_L(x) |,$ by Lemma \ref{deg}. Therefore $ | C_L(x) |= | Z(L)|+1.$ Also, $ Z(L) \subsetneq C_L(x) $ thus 
$| C_L(x) |\geq 2| Z(L)|. $ So, 
\begin{equation*}
0=| C_L(x) |- | Z(L)|-1\geq 2| Z(L)|-| Z(L)|-1=| Z(L)|-1.
\end{equation*} 
Hence $ | Z(L)|=1. $ By contrary,  let  $ q\geq 3. $ If $ x\in V(G_L), $ then  
$ A=\lbrace 2x, 3x, \dots, (q-1)x\rbrace $ is a non-empty subset of  $ V(G_L). $ Also, $ [x, \alpha x]=0 $ for all $ 1\leq \alpha \leq q-1$ and so none of the vertices of $A$ is not adjacent to $x.$ Therefore
$ \deg(x) < | V(G_L)|-1, $ which is a contradiction. Hence $ q=2. $
\end{proof}
\begin{cor}
Let $ L $ be a Lie algebra such that  $ | Z(L) |\neq 1 $  or $ q\neq 2. $ Then $ \diam(G_L)=2. $
\end{cor}
%\begin{prop}
%Let $ L $ be a Lie algebra with non-zero centre. Then $ G_L $ is not complete . 
%\end{prop}
%\begin{proof}
%Suppose that $ x\in V(G_L). $ Then there is $ y\in V(G_L) $ such that $ [x, y]\neq 0. $ Therefore $ [x+z, y]\neq 0 $ for every $z \in Z(L)$.
%Hence $ x+z \in V(G_L)$ and $ [x, x+z]=0 $ for all $ z\in Z(L). $ It implies that $ G_L $ is not complete. 
%\end{proof}
In the next lemma and propositions, we obtain some properties of the non-commuting graph $ G_L. $
\begin{lem}\label{lem28}
 Let $ L $ be a Lie algebra.  Then  $ \deg (x)\geq 2 $ for all $ x\in V(G_L). $
 % the degree of all vertices  of  $ G_L $  is bigger than $ 2. $
\end{lem}
\begin{proof}
First, we show that $ G_L $ has at least three vertices. Since $ L $ is non-abelian, then $ \dim L/Z(L)\geq 2. $ Hence $| L |/ | Z(L) |\geq q^{2}\geq 4 $ and so $ | L |\geq 4| Z(L) |. $ 
Also, we know that $ | V(G_L) |=| L |-| Z(L) | $ thus $| V(G_L) |\geq 3.  $
On the other hand, $ \deg (x)=| L |-| C_L(x) | $ for all $ x\in V(G_L) $ by Lemma \ref{deg}, $ | L |\geq 2| C_L(x) |  $ and $ C_L(x) $ has at least two elements $ 0 $ and $ x $ therefore $\deg (x)\geq 2.  $  
\end{proof}
\begin{cor}
Let $ L $ be a Lie algebra and   $ G_L $  be the non-commuting graph of $ L. $ Then $ G_L $ is not tree and star graph.
\end{cor}
\begin{prop}
Let $ L $ be a Lie algebra. Then $ G_L $ is Hamiltonian.
\end{prop}
\begin{proof}
	We know that  $ \deg(x)=| L|-| C_L(x)| $  for all  $ x \in V(G_L),$
by Lemma \ref{deg} and $ | V(G_L)|  \geq 3, $ by the proof of Lemma \ref{lem28}.
 Since $ q | C_L(x)|\leq | L| $
 for all $ x \in V(G_L)$ and $ q\geq 2, $ then 
\begin{equation*}
 2\deg(x)=2(| L|-| C_L(x)|)> | L|> | L|- | Z(L)|=| V(G_L) |. 
\end{equation*}
 Therefore $ \deg(x)>| V(G_L) |/2 $ for all $ x\in  V(G_L)$ and so $ G_L $ is Hamiltonian, by Theorem \ref{Dir}.
\end{proof}
\begin{prop}
Let $ L $ be a Lie algebra. Then $ G_L $ is an Eulerian graph.
\end{prop}
\begin{proof}
Suppose that $ x $ is an arbitrary vertex of  $ G_L,$  $ | L|=q^{n} $ and $  | C_L(x)|=q^{m}$ such that $ n>m. $ Since $ \deg(x)=| L|-| C_L(x)|,$ by Lemma \ref{deg},  
 we have  $ \deg(x)=q^{n}-q^{m}. $    Thus if $ q $ is either even or odd,  then in both cases $ \deg(x) $ is even.
Hence the degree of all vertices of $ G_L $ is even.
It implies that $ G_L $ is Eulerian, by Theorem \ref{eul}.
\end{proof}
\begin{prop}
Let $ L $ be an $ n $-dimensional Lie algebra over the field $ \mathbb{F}_q $ with the derived subalgebra of dimension $ 1. $ Then $ G_L $ is $ (q^{n}-q^{n-1}) $-regular.
\end{prop}
\begin{proof}
We claim that $ \dim C_L(x)=n-1 $ for all $ x\in V(G_L). $
 %Consider the Linear transformation $ ad_x: L\rightarrow L $ given by $ x\mapsto [x, y]. $ 
 By  Rank-Nullity Theorem, we have   $ \dim L=\dim \ker ad_x+\dim \operatorname{Im}ad_{x},$ where 
$ ad_x: L\rightarrow L $ is a linear transformation given by $ x\mapsto [x, y] $ for every $y \in L.$
Since $\ker ad_x=C_L(x)  $ and $ \dim L=n, $ then
$ \dim C_L(x)+\dim \operatorname{Im}ad_{x}=n. $ On the other hand, $ \operatorname{Im}ad_{x}\subseteq L^{2} $ and $ \dim L^{2}=1 $ thus $ \dim  \operatorname{Im}ad_{x}=1$ for all  $ x\in V(G_L). $ Therefore $ \dim C_L(x)=n-1 $ for all $ x\in V(G_L). $
Hence $ deg(x)=| L|-| C_L(x)|=q^{n}-q^{n-1} $ for all $ x\in V(G_L) $ and so $ G_L $ is $ ( q^{n}-q^{n-1})$-regular.
\end{proof} 
\begin{prop}
Let $ L $ be a Lie algebra. Then $ G_L $ is not complete bipartite.
\end{prop}
\begin{proof}
By contrary, let $ G_L $ be complete bipartite
and the vertex set of $ G_L $ be divided to two distinct parts $ A $ and $ B $ such that $ | A|\leq | B|.$ Then
\begin{equation*}
2| A|  \leq  | A|+| B|=|V(G_L)|=| L|-| Z(L)|
\end{equation*}
and so $  | A|  \leq  (| L|-| Z(L)|)/2.$
 Assume that $ x\in V(G_L) $ 
  such that $ x\in B. $ Then $ \deg(x)=| A| $ and we have $  | L|-| C_L(x)|=\deg(x)=| A|\leq  (| L|-| Z(L)|)/2.$ Therefore
  \begin{align}\label{eqbi}
| L|\leq 2| C_L(x)|-| Z(L)|.
\end{align}
   Since $ Z(L)\lneqq C_L(x), $  we may assume that $  | C_L(x) |/| Z(L)|=q^{t} $ for some $ t\geq 1. $ Hence 
$
|L |/|  C_L(x) |\leq 2- | Z(L) |/| C_L(x)|=2-\frac{1}{q^{t}}< 2
$
by \eqref{eqbi}. It implies that 
   %$ | L|=| C_L(x)|  $ 
$ L= C_L(x).  $ Hence $ x\in Z(L), $   which is a contradiction.
\end{proof}
\begin{prop}\label{lem do}
Let $ L $ be a Lie algebra over the field $ \mathbb{F}_q $ and $ G_L $ be a non-commuting graph. If $ \gamma (G_L)=1, $ then  
$ | Z(L)|=1 $ and $ q=2. $
\end{prop}
\begin{proof}
Since $ \gamma (G_L)=1, $  there exists $ x\in V(G_L) $ such that $ deg(x)=| V(G_L)|-1. $ Also, there is a vertex $ y $ such that
$ [x, y]\neq 0, $ by Lemma \ref{single}. By contrary, let $ | Z(L)|\geq 2. $ Then we can consider  a non-zero element $ z\in Z(L). $ 
On the other hand, $ [x+z, y]\neq 0 $ thus $ x+z \in V(G_L), $ but $ [x, x+z]=0. $ It is a contradiction with $ \deg(x)=| V(G_L)|-1. $ Therefore $ | Z(L)|=1. $ Now, we prove $ q=2. $ On contrary, let  $ q\geq 3. $ Since $ x $ is a vertex and $ q\geq 3, $ we have
$ A=\lbrace 2x, 3x, \dots, (q-1)x\rbrace \subseteq V(G_L).  $ Also, none of the vertices of $ A $ is not adjacent to $ x. $
So, $ \deg(x) < | V(G_L)|-1 $ and it is a contradiction. Hence $ q=2. $
\end{proof}
%In the following example, we show that the converse of Proposition \ref{lem do} does not hold.
Note that one can prove  Theorem \ref{thm2.8} in another way by the above proposition.\\
%\begin{proof}
%Since the domination   number of a complete graph is equal to $ 1 $ so,    $| Z(L)|=1 $  and  $ q=2$ by Proposition \ref{lem do}
%\end{proof}
In the following example, we show that the converse of previous proposition is not true.
\begin{exa}
Let $ L=\langle x_1, y_1, x_2, y_2 \mid [x_1, y_1]=x_1,  [x_2, y_2]=x_2 \rangle$ is a Lie algebra over the field $ \mathbb{F}_2. $ It is 
obvious $ | Z(L)|=1 $ and  $ | L|=16. $  So, $ | V(G_L)|=15 $ and $ V(G_L)=L\setminus \lbrace 0 \rbrace. $ According to the representation
of $ L $ and the vertices of $ G_L, $ we can see that
every vertex is not adjacent to at least a vertex expect itself. Hence $ \gamma (G_L)\neq 1. $
\end{exa}
\begin{thm}
Let $ L $ be a Lie algebra and $ G_L $ be the non-commuting graph  graph associated of $L$. Then $ \gamma (G_L)=1 $ if and only if there is $ x\in V(G_L) $ 
such that $ | C_L(x)|=2. $
\end{thm}
\begin{proof}
Let $ \gamma (G_L)=1. $ Then there is  $ x\in V(G_L) $ such that $ \deg(x)=| V(G_L)|-1. $ On the other hand, $ \deg(x)=| L|- | C_L(x)|, $ by Lemma \ref{deg} and
$ | V(G_L)|=| L|- | Z(L)|, $ we have $ | L|- | C_L(x)|=| L|- | Z(L)|-1.  $ So,  $  | C_L(x)|=| Z(L)|+1.  $ Also, $ | Z(L)|=1, $ by 
Proposition \ref{lem do}, hence  $ | C_L(x)|=2. $ Conversely, let $ | C_L(x)|=2$ for some $ x\in V(G_L). $ Then $ C_L(x)=\lbrace 0, x \rbrace. $
Hence $ x $ is adjacent to all vertex of $ G_L $ expect itself and $ \gamma (G_L )=1. $
\end{proof}
\section{planarity and outerplanarity of the non-commuting graph $ G_L $}
In this section, we  classify all planar and outerplanar non-commuting graphs of a Lie algebra $L.$
First, we prove the following propositions that play an important role to determine planar graphs.
\begin{lem}\label{1}
There is no  $3$-dimensional Lie algebra over the field $ \mathbb F_{2} $ with the derived subalgebra of dimension $1$ and  $Z(L) = 0.$
\end{lem}
\begin{proof}
By contrary, let $ L $ be a $3$-dimensional Lie algebra over the field $ \mathbb F_{2} $ such that $ \dim L^{2}=1 $ and $Z(L) = 0.$ Assume that $ \lbrace x, y, z \rbrace $ is a basis of $ L. $ Then the set $ \lbrace [x, y], [x, z], [y, z]\rbrace $ generates $ L^{2}. $
Without loss generality, let $\langle [x, y] \rangle $ be a basis of $ L^{2}.$
Then $ [x, z]=\alpha [x, y] $ and $ [y, z]=\beta [x, y] $ for $ \alpha, \beta \in \mathbb{F}_2. $ Since $ Z(L) = 0, $ then at least one
of the coefficients $ \alpha $ and $ \beta $ is non-zero. If $ \alpha=1 $ and $ \beta=0, $ then we can check $ y+z\in Z(L). $
It is a contradiction. Also, $ x+z\in Z(L) $ when $ \alpha=0 $ and $ \beta=1 $ and so we have a contradiction.
Finally, if $ \alpha=\beta=1, $  then $ x+y+z\in Z(L). $ Again, it is a contradiction. Hence there is no such Lie algebra.
\end{proof}
\begin{prop}\label{2}
 Let $ L $ be a $3$-dimensional Lie algebra over the field $ \mathbb F_{2} $ with $Z(L) = 0 $ and $ \dim L^{2}=2. $
Then $G_L $ is isomorphic to the following graph:
\begin{center}
\begin{figure}[H]
\centering
\begin{tikzpicture}
[scale=0.75]
\node [draw,circle,fill=black,inner sep=1pt,label=above:\tiny{$y+z$}] (C00) at (3.5,4) {};
\node [draw,circle,fill=black,inner sep=1pt,label=right:\tiny{$z$}] (C0) at (4.75,3) {};
\node [draw,circle,fill=black,inner sep=1pt,label=left:\tiny{$y$}] (C1) at (2.25,3) {};
\node [draw,circle,fill=black,inner sep=1pt,label=left:\tiny{$x$}] (C2) at (2,2) {};
\node [draw,circle,fill=black,inner sep=1pt,label=right:\tiny{$x+y$}] (C3) at (4.5,1) {};
\node [draw,circle,fill=black,inner sep=1pt,label=left:\tiny{$x+z$}] (C4) at (2.5,1) {};
\node [draw,circle,fill=black,inner sep=1pt,label=right:\tiny{$x+y+z$}] (C5) at (5,2) {};
\draw (C00)--(C2)--(C3)--(C4)--(C5)--(C00)--(C2)--(C4)--(C0)--(C2)--(C3)--(C0)--(C5)--(C3)--(C00)--(C5)--(C2);
\draw (C00)--(C4)--(C1)--(C3);
\draw (C5)--(C1)--(C2);
\end{tikzpicture}
\caption{}\label{f1}
\end{figure}
\end{center}
\end{prop}
\begin{proof}
First, let the set $ \lbrace y, z \rbrace $ be a basis of $ L^{2} $ and expand it to $ \lbrace x, y, z\rbrace $ of a basis of $ L. $
Then $ \lbrace [x, y], [x, z], [y, z] \rbrace $ generates $ L^{2}. $
Since $ L^{2} $ is abelian by \cite[Lemma 3.3]{wil}, we have $ [y, z]=0 .$ On the other hand, $ \dim L^{2}=2 $ thus 
$ \lbrace [x, y], [x, z] \rbrace$ is another basis of $ L^{2}. $ Also, we know that $ V(G_L)=\lbrace x, y, z, x+y, x+z, y+z, x+y+z \rbrace. $  Therefore  
the set  $ \lbrace [x, y], [x, z] \rbrace $ is linearly independent. So, one see that the brackets
$ [x, y],[x, z], [x, x+y], [x, x+z], [x, y+z] $ and $[x, x+y+z] $ are non-zero. Therefore $ x $ is adjacent to 
$ y, z, x+y, x+z, y+z $ and $ x+y+z. $ By a similar method, adjacency of other vertices can be found.
Hence $ G_L $ is isomorphic to Figure \ref{f1}.
\end{proof}
\begin{prop}\label{3}
 Let $ L $ be a $3$-dimensional Lie algebra over the field $ \mathbb F_{2} $ with $Z(L) = 0 $ and $ \dim L^{2}=3. $
Then $G_L $ is isomorphic to the following graph:
\begin{center}
\begin{figure}[H]
\centering
\begin{tikzpicture}[scale=0.75]
\node [draw,circle,fill=black,inner sep=1pt,label=above:\tiny{$x$}] (C00) at (3.5,4) {};
\node [draw,circle,fill=black,inner sep=1pt,label=above:\tiny{$z$}] (C0) at (4.75,3) {};
\node [draw,circle,fill=black,inner sep=1pt,label=above:\tiny{$y$}] (C1) at (2.25,3) {};
\node [draw,circle,fill=black,inner sep=1pt,label=left:\tiny{$x+y$}] (C2) at (2,2) {};
\node [draw,circle,fill=black,inner sep=1pt,label=right:\tiny{$y+z$}] (C3) at (4.5,1) {};
\node [draw,circle,fill=black,inner sep=1pt,label=left:\tiny{$x+z$}] (C4) at (2.5,1) {};
\node [draw,circle,fill=black,inner sep=1pt,label=right:\tiny{$x+y+z$}] (C5) at (5,2) {};
\draw (C00)--(C0)--(C1)--(C2)--(C3)--(C4)--(C5)--(C00)--(C2)--(C4)--(C0)--(C2)--(C3)--(C0)--(C5)--(C3)--(C00)--(C1)--(C5)--(C2);
\draw (C00)--(C4)--(C1)--(C3);
\end{tikzpicture}
\caption{}\label{f2}
\end{figure}
\end{center}
\end{prop}
\begin{proof}
Let $ \lbrace x, y, z\rbrace $ be a basis of $ L. $ Since $ Z(L)=0, $ so $V(G_L)=\lbrace  x, y, z, x+y, x+z, y+z, x+y+z \rbrace.$
Also, $ \dim L^{2}=3 $ and $ B=\lbrace [x, y], [x, z], [y, z] \rbrace $ generates 
$ L^{2}, $ then $ B $ is a basis of $ L^{2}. $ Hence all elements of $ V(G_L) $
are adjacent and $ G_L $ is isomorphic to Figure \ref{f2}. 
\end{proof}
\begin{prop}\label{4}
 Let $ L $ be a $3$-dimensional Lie algebra over the field $ \mathbb F_{2} $ with $\dim Z(L) =1.  $ Then $ \dim L^{2}=1 $ and 
$ G_L $ is isomorphic to the following graph:
\begin{center}
\begin{figure}[H]
 \begin{tikzpicture}[scale=0.75]
\node [draw,circle,fill=black,inner sep=1pt,label=above:\tiny{$x$}] (C0) at (4,4) {};
\node [draw,circle,fill=black,inner sep=1pt,label=above:\tiny{$y$}] (C1) at (3,4) {};
\node [draw,circle,fill=black,inner sep=1pt,label=left:\tiny{$x+y$}] (C2) at (2.25,3) {};
\node [draw,circle,fill=black,inner sep=1pt,label=right:\tiny{$y+z$}] (C3) at (4,2) {};
\node [draw,circle,fill=black,inner sep=1pt,label=left:\tiny{$x+z$}] (C4) at (3,2) {};
\node [draw,circle,fill=black,inner sep=1pt,label=right:\tiny{$x+y+z$}] (C5) at (4.75,3) {};
\draw (C0)--(C1)--(C4)--(C3)--(C5)--(C0)--(C2)--(C4)--(C5)--(C1)--(C2)--(C3)--(C0);
\end{tikzpicture}
\caption{}\label{f3}
\end{figure}
\end{center}
\end{prop}
\begin{proof}
Suppose that $ Z(L)=\langle z \rangle $ thus we expand $ \lbrace z \rbrace $ to a basis of $ \lbrace x, y, z \rbrace $ of $ L. $
Therefore $ \lbrace [x,y], [x,z], [y,z] \rbrace$ generates $ L^{2}. $ On the other hand, $ z\in Z(L) $ and so $ [x, z]=[y, z]=0. $
It implies $ L^{2}=\langle [x, y] \rangle $  and so $ [x, y]\neq 0. $ Also, $ V(G_L)=L\ Z(L)=\lbrace x, y, z, x+y, x+z, y+z, x+y+z \rbrace.$
 Since 
$ [x, y], [x, z], [x, x+y], [x, y+z]$ and  $[x, x+y+z] $ are non-zero, then $ x $ is adjacent to the vertices $ y, z, x+y, y+z $ and
$ x+y+z. $ By a similar way, the adjacency of other vertices are obtained .
Hence $ G_L $ is isomorphic to Figure \ref{f3}.
\end{proof}
\begin{thm}\label{5}
Let $ L $ be a $ 3 $-dimensional Lie algebra over the field $ \mathbb{F}_2. $ Then $ G_L $ is isomorphic to one of the Figures 
\ref{f1}, \ref{f2} or \ref{f3}.
\end{thm}
\begin{proof}
Since $ L $ is non-abelian, we have $ \dim L^{2}\neq 0 $ and $ \dim L/ Z(L)\geq 2. $ Also, $ \dim L=3 $ thus $ 1\leq \dim L^{2}\leq 3 $  and $ \dim Z(L)\leq 1. $
Now, the proof is completed by Lemma \ref{1}, Propositions \ref{2}, \ref{3} and \ref{4}.
\end{proof}
\begin{thm}\label{plan}
Let $ G_L $ be the non-commuting associated graph to a Lie algebra $ L. $  Then $ G_L $ is planar if and only if it is isomorphic to one of the following graphs:
\begin{figure}[H]
\begin{center}
\begin{minipage}{.4\textwidth}
\begin{tikzpicture}
 [scale=1]
\node [draw,circle,fill=black,inner sep=1pt,label=above:\tiny{$x$}] (C0) at (2,2) {};
\node [draw,circle,fill=black,inner sep=1pt,label=left:\tiny{$y$}] (C1) at (1,1) {};
\node [draw,circle,fill=black,inner sep=1pt,label=right:\tiny{$x+y$}] (C2) at (3,1) {};
\path (C0) edge  (C1);
\path (C0) edge  (C2);
\path (C1) edge  (C2);
\end{tikzpicture}
\caption{}
\label{f4}
\end{minipage}
%\hspace{cm}
\begin{minipage}{.4\textwidth}
\begin{tikzpicture}
 [scale=0.75]
\node [draw,circle,fill=black,inner sep=1pt,label=above:\tiny{$x$}] (C0) at (4,4) {};
\node [draw,circle,fill=black,inner sep=1pt,label=left:\tiny{$y$}] (C1) at (3,3) {};
\node [draw,circle,fill=black,inner sep=1pt,label=left:\tiny{$x+z$}] (C2) at (3,2) {};
\node [draw,circle,fill=black,inner sep=1pt,label=right:\tiny{$x+y+z$}] (C3) at (5,3) {};
\node [draw,circle,fill=black,inner sep=1pt,label=right:\tiny{$y+z$}] (C4) at (5,2) {};
\node [draw,circle,fill=black,inner sep=1pt,label=right:\tiny{$x+y$}] (C5) at (4,1) {};
\path (C0) edge  (C1);
\path (C1) edge  (C3);
\path (C3) edge  (C4);
\path (C4) edge  (C2);
\path (C2) edge  (C1);
\path (C4) edge  (C5);
\path (C5) edge  (C2);
\path (C0) edge  (C3);
\draw (C1) to [bend right=66]  (C5);
\draw (C5) to [bend left=45] (55:3) to [bend left=45]  (C0);
\draw (C4) to [bend right=66] (C0);
\path (C2) edge  (C3);
\end{tikzpicture}
\caption{}
\label{f5}
\end{minipage}
\end{center}
\end{figure}
\end{thm}
\begin{proof}
Suppose that $ G_L $ is a planar graph thus there is  $ x \in V(G_L) $ such that $ deg(x)\leq5, $ by Theorem \ref{deg1}. 
Since $ \deg(x)=| L|-| C_L(x)|, $ by Lemma \ref{deg} and
$ | C_L(x)|\leq | L|/q, $ then  $ | L|\leq 5q/q-1. $ On the other hand, $ f(q)=5q/q-1 $ for all $ q\geq 2 $ is a descending function 
thus $ | L|\leq 10. $ 
It is obvious that $ | L| $ is $ 2^{2}, $ $ 2^{3} $ or $ 3^{2}. $ Consider the following cases:\\
\textbf{Case 1.} $ | L|=2^{2}. $  Since $ L $ is a non-abelian Lie algebra, we have
$ \dim L/Z(L)\geq 2 $ and so $ | Z(L)|  $ can be $ 1. $ Let $ L=\langle x, y \rangle $ and $ | Z(L)| =1.$ Then 
$ V(G_L)=\lbrace x, y, x+y \rbrace, $  $ [x, y]\neq 0, $ $ [x, x+y]\neq 0 $ and $ [y, x+y]\neq 0. $ Hence $ G_L $ is isomorphic to 
 Figure \ref{f4}. \\
 %If $ | Z(L)|=2,  $ then by considering $ \lbrace x \rbrace $ as a basis for $ Z(L), $ we can expand it to the basis 
% $ \lbrace x, y \rbrace $ of $ L. $ Since $ x $ is a central element, $ L $ is an abelian Lie algebra, which is a contradiction.\\
  \textbf{Case 2.} $ | L|=3^{2}. $ By a similar method of  case 1, we can see that  $ | Z(L)|  $ is $ 1 $ or $ 3. $ Let
  $ \lbrace x, y \rbrace $ be a basis for $ L. $ Then $ L^{2} = \langle [x,y] \rangle. $ Since $ L $ is non-abelian, then $ [x,y] \neq 0. $ Then $ V(G_L)=\lbrace x, 2x, y, 2y, x+y, 2x+y, x+2y, 2x+2y \rbrace. $ 
 So, the graph $ G_L $ contains the subgraph $ K_{3,3} $ as following:
\begin{figure}[H]
\centering
\begin{tikzpicture}
[scale=0.75]
\node [draw,circle,fill=black,inner sep=1pt,label=left:\tiny{$2x+y$}] (C0) at (1,1) {};
\node [draw,circle,fill=black,inner sep=1pt,label=left:\tiny{$y$}] (C1) at (1,2) {};
\node [draw,circle,fill=black,inner sep=1pt,label=left:\tiny{$2y$}] (C2) at (1,3) {};
\node [draw,circle,fill=black,inner sep=1pt,label=right:\tiny{$ x+y$}] (C3) at (5,3) {};
\node [draw,circle,fill=black,inner sep=1pt,label=right:\tiny{$ x$}] (C4) at (5,2) {};
\node [draw,circle,fill=black,inner sep=1pt,label=right:\tiny{$ 2x+2y$}] (C5) at (5,1) {};
\draw (C2)--(C3)--(C1)--(C4)--(C0)--(C5)--(C2)--(C4)--(C0);
\draw (C1)--(C5);
\draw (C0)--(C3);
\end{tikzpicture}
\caption{}\label{f6}
\end{figure}
 Hence $ G_L $ is not planar, by Theorem \ref{Lo}. \\
 %Assume that  $ | Z(L)|=3  $  and $ \lbrace x\rbrace $ is a basis of $ Z(L). $ The set 
% $ \lbrace x, y\rbrace $ is a basis of $L$ by expanding  $ \lbrace x\rbrace. $  Since $ x $ a central element, $ L $ is an abelian Lie algebra. 
%It is a contradiction.\\
\textbf{Case 3.}  $ | L|=2^{3}. $ In this case,  $ G_L $ is isomorphic to Figures 
\ref{f1}, \ref{f2} or \ref{f3}, by Theorem \ref{5}.  
If $ G_L $ be isomorphic to Figure
\ref{f3},  then we can draw the graph by a another way and see that it is planar. So, $ G_L $ is isomorphic to Figure \ref{f5}.
Let $ G_L $ is isomorphic to Figure
\ref{f1}. Then
the graph has the following subgraph  that is isomorphic to $ K_{3,3}. $ 
\begin{figure}[H]
\centering
\begin{tikzpicture}
[scale=0.75]
\node [draw,circle,fill=black,inner sep=1pt,label=left:\tiny{$z$}] (C0) at (1,1) {};
\node [draw,circle,fill=black,inner sep=1pt,label=left:\tiny{$y$}] (C1) at (1,2) {};
\node [draw,circle,fill=black,inner sep=1pt,label=left:\tiny{$x$}] (C2) at (1,3) {};
\node [draw,circle,fill=black,inner sep=1pt,label=right:\tiny{$ x+z$}] (C3) at (5,3) {};
\node [draw,circle,fill=black,inner sep=1pt,label=right:\tiny{$ x+y$}] (C4) at (5,2) {};
\node [draw,circle,fill=black,inner sep=1pt,label=right:\tiny{$ x+y+z$}] (C5) at (5,1) {};
\draw (C2)--(C3)--(C1)--(C4)--(C0)--(C5)--(C2)--(C4)--(C0);
\draw (C1)--(C5);
\draw (C0)--(C3);
\end{tikzpicture}
\caption{}\label{f7}
\end{figure}
Hence $ G_L $ is not planar by Theorem \ref{Lo}. 
Again, $ G_L $ is not planar by Theorem \ref{Lo} when it is isomorphic to Figure \ref{f2}.    Therefore $ G_L $ is isomorphic to Figures \ref{f4} or \ref{f5}.
\end{proof}
In the next theorem, we study the existence of  non-commuting outerplanar graph associated to a Lie algebra $ L. $
\begin{thm}
Let $L$ be a Lie algebra and $G_L$ be  the non-commuting associated graph  to  $ L. $ If $G_L$ is outerplanar if and only if  $G_L$ is isomorphic to
Figure \ref{f4}.
\end{thm} 
\begin{proof}
Let $ G_L $ be an outerplanar graph. Then there is a vertex $ x $ such that $ deg(x)\leq 2 $ by Theorem \ref{outer}. On the other hand, we know that
$ deg(x)= | L|-| C_L(x)|\leq 2$  and $ | C_L(x)|\leq | L|/q. $ It implies that $ | L|\leq 2q/q-1 $ for all $ q\geq 2. $ Since $ f(q)=2q/q-1 $ for all
$ q\geq 2 $ is descending function, we have $ | L|\leq 4. $ So, $ | L| $ is $ 2, $ $ 3, $ or $ 2^{2}. $ If $ | L|=2, $ or $ 3, $ then $ L $ is abelian and it is a contradiction. Therefore $ | L|=2^{2}. $ Now, $ G_L $ is isomorphic to Figure \ref{f4} ,by a similar way in the case 1 of  the proof of Theorem \ref{plan}.
%Since $ \dim L/Z(L)\geq 2 $ for every Lie algebra $ L, $  $ | Z(L)| $ is $ 1 $ or $ 2. $
%If  $ | Z(L)|=1, $ then $ L \cong \langle x,y \vert [ x,y] = x \rangle $ by \cite[Theorem 3.1]{wil} and so $G_L$ is isomorphic the part A of Figure \ref{f4}. Let  $ | Z(L)|=2. $ Then similar to case (i) in the proof of  Theorem \ref{plan},  $ G_L $ is not outer planar.
The converse is obvious.
\end{proof}
\section{Isomorphism of the non-commuting graphs}
In this section, we are going to see that if $ G_{L_1}\cong G_{L_2} $ for two Lie algebras $ L_1 $ and $ L_2, $ then what properties
can be deduced  from this isomorphism. In the following theorem, we show that if there exists a vertex  of degree prime power number in
$ G_{L_1} $ and non-commuting graphs of  two Lie algebras $ L_1 $ and $ L_2 $ are isomorphic, then their fields have the same characteristic.
\begin{thm}
Let $ L_1 $ and $ L_2 $ be  two finite-dimensional Lie algebras over fields $ \mathbb{F}_{q_1}$ and $ \mathbb{F}_{q_2}, $ respectively. 
%$ G_{L_1} $ and $ G_{L_2} $ be non-commuting graphs of $ L_1 $ and $ L_2, $ 
Also, suppose that  $ G_{L_1} \cong G_{L_2} $  and there exists a vertex $ x\in V(G_{L_1}) $ such that $ \deg(x) $ is a prime power number. Then
$ q_1=q_2. $
\end{thm}
\begin{proof}
Since $G_{L_1} \overset{\phi}{\cong} G_{L_2}, $ 
we have $| V(G_{L_1})|=|  V(G_{L_2})|  $ and 
\begin{equation*}
 | L_1|-| C_{L_1}(x)|=\deg(x)=\deg(\varphi (x))=| L_2|-| C_{L_2}(\varphi (x))|,
\end{equation*}
by Lemma \ref{deg}.
Let  $ | L_1|=q_1^{n_1}, $ $ | L_2|=q_2^{n_2}, $ 
$ | C_{L_1}(x)|=q_1^{m_1} $ and   $ | C_{L_2}(x)|=q_2^{m_2}, $ where $ n_1> m_1 $ and $ n_2> m_2. $
Then
%$ q_1^{n_1}- q_1^{m_1}=q_1^{n_2}- q_1^{m_2}$ or 
$ q_1^{m_1}(q_1^{n_1-m_1}- 1)=q_2^{m_2}(q_2^{n_2-m_2}- 1).$
On the other hand, $ \deg(x)=p^{n} $ for some prime $ p $ and $ n\geq 1 $ thus $ q_1\vert p^{n} $ and $ q_2\vert p^{n}, $ which  imply that $ q_1=p=q_2 $ and the proof follows.
\end{proof}

\begin{thm}
Let $ L_1 $ and $ L_2 $ be Lie algebras and  $ G_{L_1} \cong G_{L_2}. $ If 
there exist a vertex $ x\in V(G_{L_1}) $ such that $ \deg(x)=q, $ where $ q $ is a prime number,  then $ | L_1|= | L_2|.$
\end{thm}
\begin{proof}
We know that  $ L_1 $ and $ L_2 $ are the Lie algebras over the $ \mathbb{F}_p, $ by the previous theorem.
Since $ G_{L_1} \overset{\phi}{\cong} G_{L_2}, $ 
%then there is  an isomorphism $ \varphi $ from $ V(G_{L_1}) $ to $ V(G_{L_2}). $ Thus  
we have $ q=\deg(x)=\deg(\varphi(x)) $ and so 
%$  | L_1|-| C_{L_1}(x)|=| L_2|-| C_{L_2}(\varphi (x))|=q$ or
 $| C_{L_1}(x)|(| L_1|/ | C_{L_1}(x)|-1)=| C_{L_2}(\varphi(x))|(| L_2|/ | C_{L_2}(\varphi(x))|-1)=q.  $ One can see that  $ | C_{L_1}(x)|\geq 2 $ and $ | C_{L_2}(\varphi(x))|\geq 2, $ so $ | C_{L_1}(x)|= | C_{L_2}(\varphi(x))|=q. $
On the other hand, $ Z(L_1)\lneqq  C_{L_1}(x)$ and $ Z(L_2)\lneqq  C_{L_2}(\varphi(x))$ imply that 
$  | Z(L_1)|=| Z(L_2)|=1.$ Hence $ | V(G_{L_1})|=| V(G_{L_2})| $ by the existence of isomorphism between $ G_{L_1} $ and $ G_{L_2}$  and $ | L_1|-| Z(L_1)|=| L_2|-| Z(L_2)|. $
Therefore $ | L_1|= | L_2|,$ as required.
\end{proof}

\begin{thm}
Let $ L_1 $ and $ L_2 $ be Lie algebras and $ G_{L_1} \cong G_{L_2}. $ If  
 there exist a vertex $ x\in V(G_{L_1}) $ such that $ \deg(x)=pq, $ where $ p $ and $ q $
are distinct prime numbers and $ p> q, $ then  $  L_1\cong L_2. $
\end{thm}
\begin{proof}
By a similar way of the previous theorem, 
%there is an isomorphism $ \varphi $ from  $ V(G_{L_1}) $ to $V(G_{L_2})$ and
 we have
  \begin{equation*}
  | C_{L_1}(x)|(| L_1|/ | C_{L_1}(x)|-1)=| C_{L_2}(\varphi(x))|(| L_2|/ | C_{L_2}(\varphi(x))|-1)=pq
  \end{equation*}
  and so there are three possibilities:
\\
\textbf{(i)} $ | C_{L_1}(x)|=p $ and $ | L_1|/ | C_{L_1}(x)|-1=q. $ In this case,  $ | L_1|=pq+p=(q+1)p. $ 
%On the other hand, 
%$ | L_1|=q^{n_1}, $ where $ q_1 $ is a prime number and $ n_1\geq 1. $ 
Thus the  only possibility
is when $ p=3 $ and $ q=2. $ Therefore $| L_1|=3^{2}.  $
\\
\textbf{(ii)} $ | C_{L_1}(x)|=q $ and $ | L_1|/ | C_{L_1}(x)|-1=p. $ Then $ | L_1|=pq+q=(p+1)q $ and so   $ p=2 $ and $ q=3. $
Since $ p> q, $ then it  is impossible. 
\\
\textbf{(iii)} $ | C_{L_1}(x)|=pq $ and $ | L_1|/ | C_{L_1}(x)|-1=1. $ Thus $ | L_1|=2pq, $ which is not possible. \\
Hence the only possibility is $ | C_{L_1}(x)|=p. $ Similar to the proof of previous theorem,   $ | C_{L_2}(\varphi(x))|=p, $
where $ \varphi $ is an isomorphism from $  C_{L_1}$ to $ C_{L_2}. $ Thus  $ | Z(L_1)|=| Z(L_2)|=1, $ which implies that 
$ | L_1|=| L_2|=3^{2}. $ Since up to isomorphism there is a unique $ 2 $-dimensional non-abelian Lie algebra, then
$ L_1\cong L_2. $
\end{proof}
\begin{thm}
 Let $ L_1 $ and $ L_2 $ be Lie algebras and $ G_{L_1} \cong G_{L_2}. $
If there exists a vertex $ x\in V(G_1) $ such that $ \deg(x)=p^{2}q, $ where $ p $ and $ q $
are prime numbers and $ p> q. $ Then   $  L_1\cong L_2. $
\end{thm}
\begin{proof}
We know that $ p^{n}q=\deg(x)= | L_1|-| C_{L_1}(x)|=| C_{L_1}(x)|(| L_1|/ | C_{L_1}(x)|-1).$ Thus there are the following
possibilities:\\
\textbf{(i)} $ | C_{L_1}(x)|=p. $ We have $ | L_1|/ | C_{L_1}(x)|-1=pq $ and so $ | L_1|=p^{2}q+p=p(pq+1). $ Since $| L_1|$ must be a prime power number, this is impossible.
\\
\textbf{(ii)} $| C_{L_1}(x)|=q.$ Then $ | L_1|=p^{2}q+p=q(p^{2}+1) $ and again it is impossible.
\\
\textbf{(iii)} $| C_{L_1}(x)|=pq.$ Since $| C_{L_1}(x)|$ is not a prime power, this case does not occur.
\\
\textbf{(iv)} $| C_{L_1}(x)|=p^{2}.$ We have $ | L_1|=p^{2}q+p=p^{2}(q+1) $ and the only possibility is when $ p=3 $ and $ q=2. $
Then $ | L_1|=3^{2}. $ Now, if $ \varphi $ is an isomorphism between $ G_{L_1} $ and $ G_{L_2}, $ then  
$ p^{2}q=\deg(x)=\deg(\varphi(x))=| L_2|-| C_{L_2}(\varphi(x))|=| C_{L_2}(\varphi(x))|(| L_2|/| C_{L_2}(\varphi(x))|-1). $ 
So,$ | C_{L_2}(\varphi(x))|=p^{2} $ and  $ | L_2|=3^{2}, $ 
by considering the above cases. Hence $ | L_1|=| L_2|=3^{2}. $ Since, there is only a $ 2 $-dimensional non-abelian Lie algebra, by  \cite[Theorem 3.1]{wil}, then $ L_1\cong L_2. $ 
%Now, when  $ | C_{L_1}(x)|=| C_{L_2}(\varphi(x))|=3, $ then  $ | Z(L_1)|= | Z(L_2)|=1$ or $ 3. $ 
% Therefore we have two cases:\\
%Case 1  ???
\end{proof}
 
 \begin{cor}
 Let  Let $ L_1 $ and $ L_2 $ be Lie algebras and $ G_{L_1} \cong G_{L_2}. $
 If there is a vertex $ x\in V(G_{L_1}) $ such that $ \deg(x)=p^{n}q, $ where $ p $ and $ q $ are prime numbers, $ p> q $ and $ n\geq 1. $ Then   $| L_1|=|L_2|. $ 
 \end{cor}
 \begin{proof}
 The proof is similar to the above theorems.
 \end{proof}
\begin{prop}
 Let $ L_1 $  and $ L_2 $ be    Lie algebras over the field
   $ \mathbb{F}_{2}  $ such that $ G_{L_1} \cong G_{L_2}. $ Then $| L_1|=|L_2|. $ 
\end{prop}
\begin{proof}
We know that $ | V(G_{L_1})| = | V(G_{L_2})|$ when $ G_{L_1} \cong G_{L_2}. $ Then
  $ | L_1|-|Z(L_1)|= | L_2|-|Z(L_2)|.$ 
Suppose that $ | L_1|=2^{n_1}, $ $ | L_2|=2^{n_2}, $ $ | Z(L_1)|=2^{m_1} $ and $ | Z(L_2)|=2^{m_2}. $ Thus
$ 2^{n_1}-2^{m_1}=2^{n_2}-2^{m_2}. $ We can assume that $ m_1\geq m_2. $ So, 
\begin{equation*}
2^{m_1-m_2}(2^{n_1}-1)=2^{n_2}-1.
\end{equation*}
On the other hand, $ 2^{n_2}-1 $ is odd and $ 2^{m_1-m_2}(2^{n_1}-1) $ is even thus $ 2^{m_1-m_2}=0. $ Hence $ m_1=m_2 $ and so $ | Z(L_1)| = |Z(L_2)|. $ It implies that $ | L_1| = |L_2|. $
\end{proof}
It is a common question to ask what properties of Lie algebras are preserved under isomorphism between the non-communing graphs of them. 
The next example shows that it is not true for  the property of nilpotency. 
%preserved under isomorphism between non-communing graphs of Lie algebras. 
\begin{exa}
Assume that $ L_1=\langle x, y, z\vert [x, y]=z \rangle $ and $ L_2=\langle x, y, z\vert [x, y]=x \rangle. $ We cam see that
the non-commuting graph associated to $ L_1 $ and $ L_2 $ is isomorphic to Figure \ref{f3}. Also, $ L_1 $ is nilpotent, but $ L_2 $
does not have such property. So, the nilpotency property is not kept under isomorphism of the non-communing graphs.
\end{exa}
It is interesting to see that under what sufficient and necessary conditions we have $ G_{L_1}\cong G_{L_2} $
if and only if $ | L_1| = |L_2|. $ We left it as a conjecture here.\\
\\
\textbf{Conjecture} Find necessary and sufficient conditions  such that $ G_{L_1}\cong  G_{L_2}$  
if and only if $ | L_1| = |L_2|. $ 
 
\end{document}